\documentclass[conference,10pt]{IEEEtran_ICRATmod}

\IEEEoverridecommandlockouts
\usepackage[numbers]{natbib}
\usepackage{subfig}
\usepackage{graphics} 
\usepackage{epsfig} 
\usepackage{times} 
\usepackage{amsmath} 
\usepackage{amssymb}  
\usepackage{amsthm}
\usepackage{amsfonts}
\usepackage{mathtools}
\usepackage{multirow} 
\usepackage{rotating} 
\usepackage{algorithm}
\usepackage{algorithmic}
\usepackage{graphicx}
\usepackage{float}
\usepackage{flushend}
\usepackage[table,xcdraw]{xcolor}
\usepackage{csquotes}
\usepackage{bm}
\usepackage{relsize}
\usepackage{epspdfconversion}
\usepackage{colortbl}
\usepackage{lipsum}
\usepackage[mathscr]{eucal}
\usepackage{fancyhdr}

\usepackage{tikz}

\newcommand{\Cross}{\mathbin{\tikz [x=1.4ex,y=1.4ex,line width=.2ex] \draw (0,0) -- (1,1) (0,1) -- (1,0);}}%
\usepackage{algorithm}
\usepackage[algo2e]{algorithm2e}

\usepackage{bbm}

\usepackage{hyperref}

\newcommand{\norm}[1]{\left\lVert#1\right\rVert}

\DeclareMathAlphabet\mathbfcal{OMS}{cmsy}{b}{n}

\def\compactify{\itemsep=0pt \topsep=0pt \partopsep=0pt \parsep=0pt}
\newcommand{\compress}{\itemsep=0pt \topsep=0pt \partopsep=0pt \parsep=0pt \leftmargin=30pt \labelwidth=30pt}

\let\latexusecounter=\usecounter

\begin{document}
\title{Distributionally Robust Ground Delay Programs with Learning-Driven Airport Capacity Predictions}

\author{\scriptsize
	   	   \IEEEauthorblockN{Haochen Wu, Max Z. Li}
	   	   \IEEEauthorblockA{University of Michigan\\Ann Arbor, MI, USA\\ \{\href{mailto:haocwu@umich.edu}{haocwu}, \href{mailto:maxzli@umich.edu}{maxzli}\}@umich.edu}
	   \and
	   	   \IEEEauthorblockN{Xinting Zhu, Shuchang Li, Ying Zhou, Lishuai Li}
	   	   \IEEEauthorblockA{City University of Hong Kong\\Kowloon Tong, Hong Kong\\ \{\href{mailto:xt.zhu@my.cityu.edu.hk}{xt.zhu}, \href{mailto:shuchanli2-c@my.cityu.edu.hk}{shuchanli2-c}, \href{mailto:y.zhou@my.cityu.edu.hk}{y.zhou}, \href{mailto:lishuai.li@cityu.edu.hk}{lishuai.li}\}@cityu.edu.hk}
}

\vspace{-2cm}
\maketitle

\renewcommand{\headrulewidth}{0pt}
\lhead{\textcolor{black}{ICRAT 2024}}
\rhead{\textcolor{black}{Nanyang Technological University, Singapore}}
\cfoot{\textcolor{black}{\thepage}}
\lfoot{\textcolor{white}{\thepage}}
\thispagestyle{fancy}
\pagestyle{fancy}

\fancypagestyle{firststyle}
{
    \fancyhf{}
    \lhead{\textcolor{black}{ICRAT 2024}}
\rhead{\textcolor{black}{Nanyang Technological University, Singapore}}
    \cfoot{\textcolor{black}{\thepage}}
    \lfoot{\textcolor{white}{\thepage}\\ \vspace{0.1cm} \scriptsize \parbox{0.465\textwidth}{\hrule ~~\\H. Wu was partially supported by a Departmental Fellowship from the University of Michigan. X. Zhu, S. Li, and Y. Zhou were partially supported by Hong Kong Innovation and Technology Commission Innovation and Technology Fund (GHP/145/20) and City University of Hong Kong Internal Fund (PJ9678283). }}
}

\noindent
\begin{abstract}
Strategic Traffic Management Initiatives (TMIs) such as Ground Delay Programs (GDPs) play a crucial role in mitigating operational costs associated with demand-capacity imbalances. However, GDPs can only be planned (e.g., duration, delay assignments) with confidence if the future capacities at constrained resources (i.e., airports) are predictable. In reality, such future capacities are uncertain, and predictive models may provide forecasts that are vulnerable to errors and distribution shifts. Motivated by the goal of planning optimal GDPs that are \emph{distributionally robust} against airport capacity prediction errors, we study a fully integrated learning-driven optimization framework. We design a deep learning-based prediction model capable of forecasting arrival and departure capacity distributions across a network of airports. We then integrate the forecasts into a distributionally robust formulation of the multi-airport ground holding problem (\textsc{dr-MAGHP}). We show how \textsc{dr-MAGHP} can outperform stochastic optimization when distribution shifts occur, and conclude with future research directions to improve both the learning and optimization stages.




\end{abstract}

\begin{small}{{\bfseries\itshape Keywords---Air traffic management; Ground Delay Programs (GDPs); Airport capacity prediction; Distributionally robust optimization }}\end{small}

\thispagestyle{firststyle}
\section{Introduction}  \label{sec:intro}
Congestion in air transportation systems results from demand-capacity imbalances, often stemming from airport capacity reductions. Within the US National Airspace System (NAS), the strategic implementation of Traffic Management Initiatives (TMIs) seeks to reduce the operational costs of such imbalances. A prominent example of TMIs is Ground Delay Programs (GDPs), which aims to strategically delay flights on the ground at the origin airport and mitigate costly airborne delays. The optimal implementation of GDPs is the purview of airport ground holding optimization problems, or GHPs. 

Airport capacities at future time periods play a significant role in GDP implementations. If such capacities are known, the optimal delay allocation decisions can be found by solving GHPs \cite{GHP1}. However, in practice, it is extremely difficult for traffic management decision-makers to ascertain future airport capacities (i.e., Airport Arrival Rates and Airport Departure Rates) due to myriad uncertainties. Such uncertainties stem from environmental (e.g., convective weather forecasts \cite{kicinger2016airport}) and operational (e.g., runway availability and traffic volume \cite{tittle2013airport}) sources. Although a large variety of prior work focuses on, e.g., airport runway configuration prediction \cite{avery2015predicting}, weather impact predictions \cite{kicinger2016airport}, and TMI implementation predictions \cite{vlachou2019simultaneous}, all predictions result in a probability distribution over potential outcomes. If such predictions were to be incorrect (e.g., due to distribution shifts or misspecification), the resultant GDP implementation may be sub-optimal. 



\subsection{Motivation and research problem}

In this work, we are motivated by the goal of leveraging advancements in machine learning (ML)-based predictions for airport capacities while adopting a \emph{cautiously optimistic} approach towards prescribing GDP solutions. The former acknowledges the adoption of ML models which enables data-driven predictive models for aviation (a key tenant in the Federal Aviation Administration's Information-Centric NAS vision \cite{Info_centric_NAS}). Moreover, the latter acknowledges the need to make \emph{robust} decisions, understanding that such probabilistic predictions may be incorrect.


ML has continued to see significant developments and applications in the field of aviation, particularly concerning airport and airspace operations. By using historical and real-time data, ML models can help airlines and airport operators anticipate and mitigate flight delays. ML is aiding in the development of more advanced ATM systems, which can predict traffic flows and flight delays, optimize flight paths, and improve overall airspace efficiency. The European Union's SESAR (Single European Sky ATM Research) and the United States' NextGen programs have been working on integrating ML to enhance ATM. However, most existing ML methods can only generate point estimations for prediction, which do not fully characterize the uncertainties brought by the dynamic and complex nature of air transportation systems. 

Existing methods for optimizing GDPs include the deterministic multi-airport ground holding problem (MAGHP) and the stochastic MAGHP; the former assumes deterministic airport capacities, while the latter relaxes this assumption. The feasibility of implementing ML models in real-world settings is compromised when the downstream learning-driven optimization model is highly sensitive to predictions from the ML model. This sensitivity leads to decision-making based on potentially inaccurate information (e.g., inaccurate airport capacity predictions). Consequently, it is critical for GDP solutions to acknowledge and address uncertainties in airport capacity prediction models. 


To account for distributional uncertainties in airport capacity predictions, we investigate a \emph{distributionally robust} version of the MAGHP, termed the distributionally robust multi-airport ground holding problem (\textsc{dr-MAGHP}). Distributionally robust optimization (DRO) aims to identify the optimal solution by considering the worst-case (with respect to the objective function) distribution within a predefined set of distributions. This predefined set is known as the \emph{ambiguity set}, and is constructed based on the predicted airport capacity distribution from the upstream ML model. In Section \ref{sec:methodology} we will rigorously restate relevant terms above. 


\subsection{Background and prior works}
The deterministic MAGHP, stochastic MAGHP, and other MAGHP variations have been well studied \cite{GHP1}. Stochastic implementations include using two-stage stochastic programming and chance-constrained programming \cite{GHP2,GHP4}. Moreover, \cite{GHP5} develops a data-driven control framework that not only minimizes the total cost but redistributes delays spatially, potentially improving operational recovery capabilities. \cite{GHP6} and \cite{GHP7} also take equity, fairness, and passenger-centric considerations into account. 

In the operations research literature, DRO is emerging as a rigorous way to combine elements of stochastic and robust optimization approaches. \cite{DR1} proposes moment-based DRO, constructing the ambiguity set using the first and second moments of the distribution. Alternatively, \cite{DR2} builds a Wasserstein distance-based ambiguity set: This approach considers all distributions within a Wassertein ball of radius $\epsilon$, and proposes a convex reformulation technique to recover a tractable form of the original DRO problem. \cite{DR3} studies the Wassertein ambiguity set-based distributionally robust mixed-integer program and solves it using dual decomposition methods. 

Airport capacity can be defined as the maximum sustainable throughput for arriving (the Airport Arrival Rate, or AAR) and departing (the Airport Departure Rate, or ADR) flights \cite{kicinger2012airport}. 
In contrast to declared capacities obtained from theoretical analyses or statistical approaches \cite{gilbo1993airport}, real-time capacity is dynamic and challenging to predict in advance. 
Real-time airport capacities are influenced by several interconnected operational and environmental factors \cite{choi2021artificial}. As the prediction time horizon increases, forecast uncertainty grows as well, rendering accurate long-term predictions difficult \cite{tien2018using}. A sampling of previous work include analytical approaches such as the Integrated Airport Capacity Model (IACM) \cite{kicinger2012airport}, and more recently, data-driven approaches such as \cite{cox2016probabilistic}. Other related work focused on predictions of, e.g., traffic flow \cite{murca2018identification}. Few works have been conducted on the prediction of real-time airport capacities \cite{wang2021prediction}, and nascent work has been done on distribution predictions for flight delays \cite{Wang2022}.

\section{Contributions}

In this work, our contributions are as follows:

\begin{enumerate}
    \item We develop an ML framework for providing distributional predictions of the airport capacity, rather than a single value. 
    \item We formulate and solve the \textsc{dr-MAGHP} given input capacity distributions from the ML model. 
    \item We conduct a sensitivity analysis to evaluate the performance of \textsc{dr-MAGHP} under varying levels of airport capacity overestimation. 
\end{enumerate}


\section{Methodology}  \label{sec:methodology}

\begin{figure*}[htbp]
    \centering
    \includegraphics[width=\textwidth]{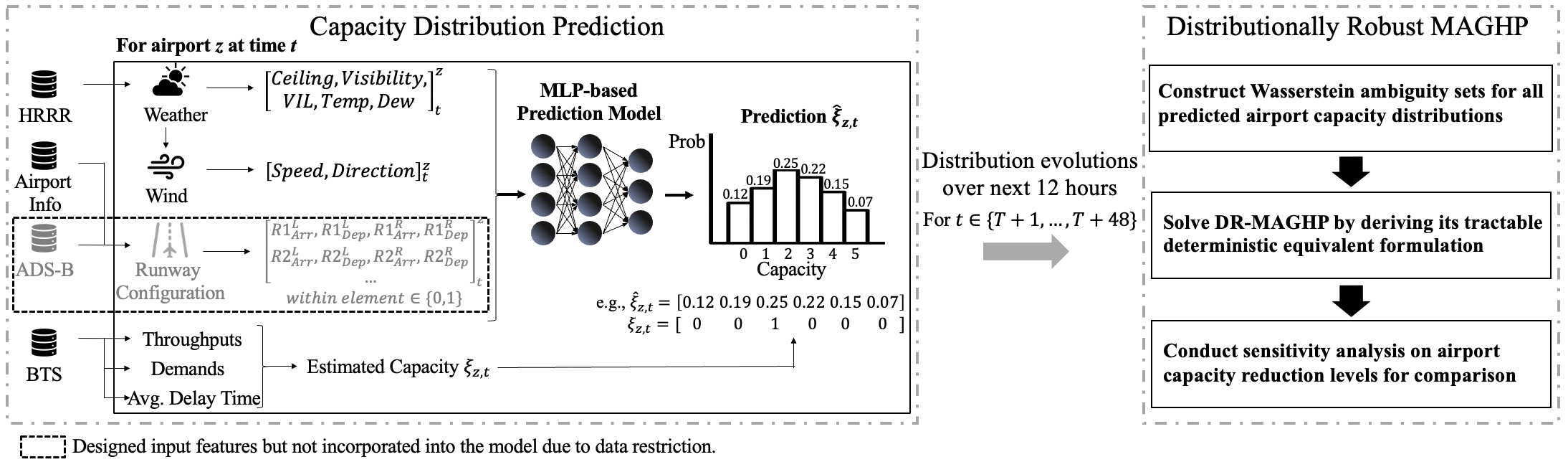}
    \caption{Learning-driven airport capacity distribution prediction and distributionally robust GDP optimization framework. 
    }
    \label{fig:capacity_distribution_prediction_framework}
\end{figure*}

We first introduce our airport capacity distribution prediction model. The predicted distributions serve as inputs to the \textsc{dr-MAGHP}. 
Figure~\ref{fig:capacity_distribution_prediction_framework} provides an overview of our combined learning-driven optimization framework. 






\subsection{Airport capacity distribution prediction}   \label{ssec:airport_cap_dis_pred}

We develop a deep learning model for airport capacity distribution prediction. The model can forecast arrival and departure capacity distributions across a 12-hour prediction horizon, discretized in 15-minute intervals. A separate prediction model is built for each airport to account for airport-specific differences \cite{liu2019using}. Furthermore, for each airport, separate and independent models are built for arrival and departure capacity predictions. 


\subsubsection{Deriving actual capacities from throughput}  \label{sssec:derive_actual_capacities}
In order to train and validate the prediction model, we need observations of the actual, true airport capacity. A unique challenge for this prediction problem is that we do not directly observe the actual capacity value from historical data \cite{ramanujam2009estimation}. 
In this study, we derive estimates of the actual capacity from historical records of airport arrival and departure throughput. The throughput at time $t$ (observable) typically are underestimations of actual capacity at time $t$ (unobservable, to be estimated), particularly during off-peak periods. To address this, we use two rules as defined in \eqref{eq:capacity_estimation} to filter out time periods with low-volume of flight operations and select peak operational times when throughput generally reaches its maximum to estimate the actual airport capacity:
\begin{equation}
\label{eq:capacity_estimation}
    \widehat{\text{capacity}}_t = \text{throughput}_t \iff \text{(Rule 1)} \lor \text{(Rule 2)},
\end{equation}
where Rule 1 is given by $\text{demand}_t \geq \text{throughput}_t + 3$ and Rule 2 by $\left(\text{avg. delay}_t > 30\right) \land \left(\text{no. of delayed flights}_t > 1 \right)$.
Only data from time periods satisfying \eqref{eq:capacity_estimation} are used to train and validate the airport capacity prediction model. 

\subsubsection{Probabilistic capacity prediction}
Our model includes input feature engineering, a multilayer perceptron (MLP)-based prediction model, and output feature engineering. 
Denote by $\widehat{\xi}_{ \{a, g\}, t }^z$ the predicted airport arrival ($a$) or departure ($g$) capacity distribution at time $t$ for airport $z$. We seek to learn a mapping $\Psi_{\{a, g\}}^z$ from inputs $X_t^z$ to capacity outputs $\widehat{\xi}_{ \{a, g\}, t }^z$, or explicitly, $\widehat{\xi}_{ \{a, g\}, t }^z = \Psi_{\{a, g\}}^z \left(X_t^z\right)$.
On the model output, the airport capacity $\widehat{\xi}_{ \{a, g\}, t }^z$, we introduce a feature engineering technique for distributional prediction: $\widehat{\xi}_{ \{a, g\}, t }^z$ is represented as a set of discrete value, encoded via a one-hot encoding method. This process converts scalar time-indexed capacity values into high-dimensional binary vector representations. For example, as shown in Figure \ref{fig:capacity_distribution_prediction_framework}, a capacity of $2$ is encoded into a sparse vector as $ {[ 0\ 0\ 1\ 0\ 0\ 0 ]}_t $, where the length of this vector is $6$, corresponding to the range of capacity values observed at this airport historically. 


Due to data limitations (e.g., historical runway configurations could not be obtained), we only use meteorological conditions for the model input $X_t$. 
Weather data are converted and vectorized, with the following seven features: ceiling, visibility, vertically integrated liquid (VIL), temperature, dew points, wind direction, and surface wind speeds. 
The selection of these features is based on their studied impacts on airport capacity (e.g., \cite{allan2001analysis,renhe2014meteorological}).

For the prediction models $\Psi_{\{a, g\}}^z$, we use a 3-layer MLP whose architecture details are given in Table~\ref{tab:MLP_setting}. The output layer of the MLP reflects the range of airport capacity classes, similar to the multiclass classification problem setting.

\begin{table}[htp]
    \centering
    \resizebox{0.45\textwidth}{!}{
    \begin{tabular}{c|c|c|c|c}
        \hline 
        \textbf{No. of layers} & \textbf{No. of neurons} & \textbf{Activation} & \textbf{Output} & \textbf{Loss} \\
        \hline
         3 & [17, 32, ${MAX}_i + 1$] & ReLU & Softmax & Cross-entropy \\
        \hline 
    \end{tabular}}
    \caption{3-layer MLP architecture for airport capacity distribution prediction. ${MAX}_i$ indicates the historical maximum capacity of airport $i$, with indexing normalized.}
    \label{tab:MLP_setting}
\end{table}

\subsection{Distributionally robust MAGHP (\textsc{dr-MAGHP})}

\begin{figure}[]
\centering
\subfloat[Gaussian]{\includegraphics[width=4.4cm]{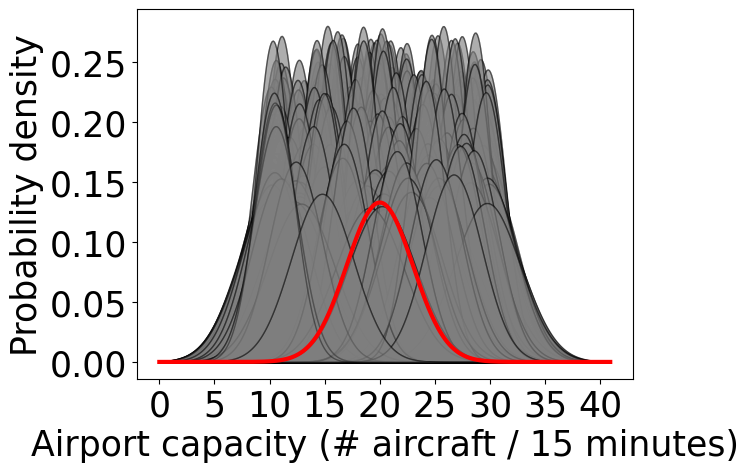}}\hfill
\subfloat[Gamma]{\includegraphics[width=4.4cm]{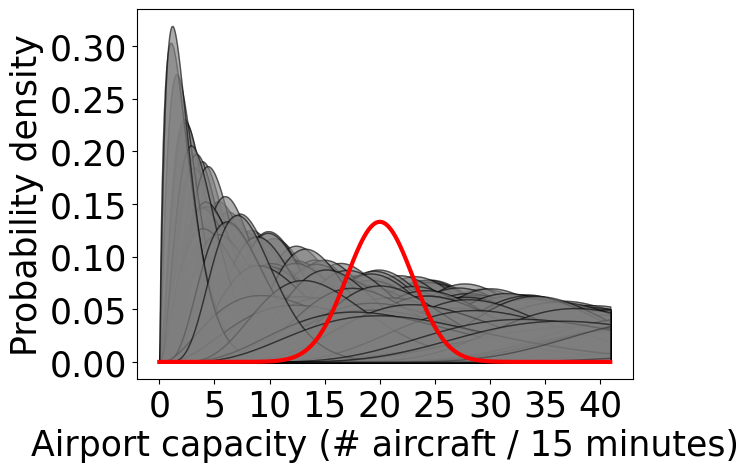}}\vfill
\subfloat[Erlang]{\includegraphics[width=4.4cm]{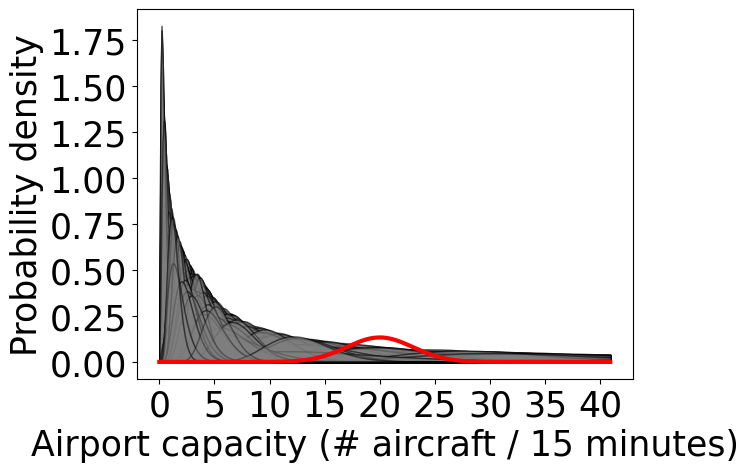}}\hfill
\subfloat[Lognorm]{\includegraphics[width=4.4cm]{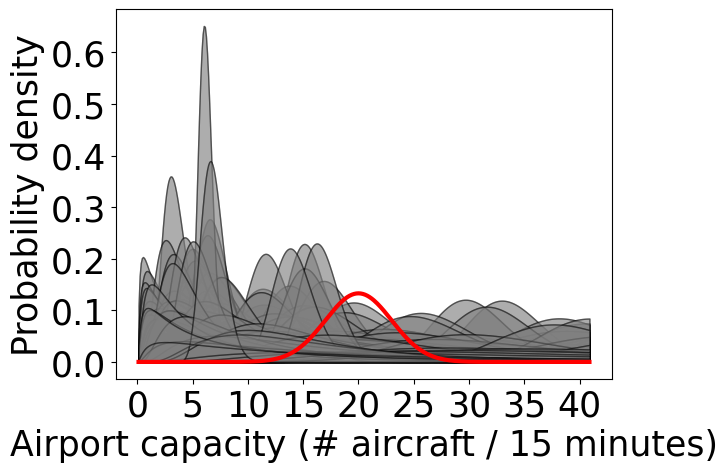}}\
\caption{Distributions in the Wasserstein ambiguity set.}\label{fig:wsdistributions}
\end{figure}

Recall that we seek to make optimal ground delay allocation decisions that are \emph{robust} with respect to the predicted airport capacity distribution (output of Section \ref{ssec:airport_cap_dis_pred}). To do so, we will solve the MAGHP under a \enquote{worst-case} capacity distribution, located within the Wasserstein ambiguity set centered at the predicted distribution. 
Figure \ref{fig:wsdistributions} shows a numerical example of different distributions---and distribution families---which lie within the ambiguity set with radius $\epsilon = 0.005$ centered at a Gaussian distribution $\mathcal{N}(20, 3)$. 
Formally, let $M(\Xi)$ be the space of all probability distributions $\mathbb{Q}$ with support $\Xi$. The Wasserstein distance $d_w : M(\Xi) \times M(\Xi) \rightarrow \mathbb{R}_{\geq 0}$ is the minimum transportation cost between distributions $\mathbb{Q}_1 \in M(\Xi)$ and $\mathbb{Q}_2 \in M(\Xi)$ \cite{DR2}, and is given explicitly as: 

\vspace{-0.3cm}
\begin{equation}
\begin{aligned}
d_{w}\left(\mathbb{Q}_1,\mathbb{Q}_2\right) = \inf_{\Pi \in \mathcal{D}_{\Pi} \left(\xi_1,\xi_2 \right)}\int_{\Xi^2} \norm{\xi_1 - \xi_2}_2 \, \Pi(d\xi_1,d\xi_2),
\end{aligned}
\label{1.0}
\end{equation}

\noindent
where $\Pi$ is a joint distribution of random variables $\xi_1$ and $\xi_2$ with marginals $\mathbb{Q}_1$ and $\mathbb{Q}_2$, respectively. 
We denote $\mathcal{D}_{\Pi}\left(\xi_1,\xi_2\right)$ as the set of all joint distributions on $\xi_1$ and $\xi_2$ with marginals $\mathbb{Q}_1$ and $\mathbb{Q}_2$. Let $Z$ be the set of all airports and $\left\{ \widehat{\xi}_1,\widehat{\xi}_2, \hdots ,\widehat{\xi}_{|Z|} \right\}$ be the set of $|Z|$ airport capacities with the corresponding estimated probabilities of occurrence $\left\{ \widehat{p}_1,\widehat{p}_2, \hdots ,\widehat{p}_{|Z|} \right\}$ (i.e., this is precisely the predicted capacity distribution's probability mass function, or PMF). The ambiguity set centered around a predicted capacity distribution $\widehat{P}$ with radius $\epsilon > 0$, denoted as $\mathcal{P}_{\epsilon}\left( \widehat{P} \right)$, is given by
\begin{equation}
\begin{aligned}
\mathcal{P}_\epsilon\left(\widehat{P}\right) \coloneqq \left\{ \mathbb{Q} \in M(\Xi) : d_{w}\left(\widehat{P},\mathbb{Q} \right) \leq \epsilon  \right\}.
\end{aligned}
\label{1.1}
\end{equation}
Note that to differentiate between arrival and departure capacities, we will use subscripts and superscripts with $a$ and $g$, respectively, when we write out the full \textsc{dr-MAGHP} model. We refer readers to \cite{drGHP} for additional technical details on specifying the ambiguity set.

We denote the set of flights by $F$, 
$D_z(t)$ and $R_z(t)$ as the departure and arrival capacity of airport $z \in Z$ at time $t$, respectively, and $\mathcal{C} = F \times F$ as the set of all flight pairs connected by the same aircraft (or tail), where $(f, f') \in \mathcal{C}$ denotes the preceding flight $f$ and successive flight $f'$. 
The decision variables $u_{f,t}$ and $v_{f,t}$ are binary variables where $u_{f,t}$ equals one if flight $f$ will departure at time $t$; similarly, $v_{f,t}$ equals one if flight $f$ will land at time $t$. $z^{g}_{f}$ and $z^{a}_{f}$ represent flights scheduled to depart or land at airport $z$, respectively. $d_f$ is the scheduled departure time of $f$ and $r_f$ is the scheduled arrival time of $f$. $T^{d}_{f}$ is the set of available time periods for $f$ to take off and $T^{r}_{f}$ is the set of available time periods for $f$ to land. $g_{f},a_{f}$ are ground holding delay and airborne delay, respectively, where $g_{f} = \sum_{t\in T^{d}_{f}}tu_{f,t} - d_f$ and $a_{f} = \sum_{t \in T^{a}_{f}}tv_{f,t} - r_{f} - g_{f}$. These definitions follow the standard formulation of the MAGHP \cite{GHP1}.

We incorporate predicted airport departure and arrival capacity distributions from Section \ref{ssec:airport_cap_dis_pred}, along with robustness guarantees, through a two-stage formulation. 
In the two-stage setting, in addition to the first stage decision variables $u_{f,t}$ and $v_{f,t}$, we introduce second stage decision variables $y^{(g)}_{z,t}$ and $y^{(a)}_{z,t}$.
Decision variables \(y^{(g)}_{z,t}\) denote the additional number of flights entering the departure queue at time period \(t\), while \(y^{(a)}_{z,t}\) represents the additional number of flights entering the arrival queue at the same time period. The second stage of the model focuses on reallocating the actual departure and arrival times of flights once the actual airport capacity distributions are realized. The objective of the second stage model is to minimize the number of flights joining either queues. With this, the full \textsc{dr-MAGHP} can be written as follows:
\begin{subequations}
\begin{alignat}{2}
\min_{u,v} \quad \left\{\sum_{f \in F} \left(C_g g_f + C_a a_f \right) + \right.&\left.  \max_{p \in \mathcal{P}_{\epsilon_g}  \left(\widehat{P}^{\left(g\right)}\right)} \, \mathbb{E}_{p}\left[Q \left(u,\xi^{(g)}\right)\right] \notag \right. \\ \left.  + \max_{p \in \mathcal{P}_{\epsilon_a}\left(\widehat{P}^{\left(a\right)}\right)} \, \mathbb{E}_{p}\left[Q\left(v,\xi^{(a)} \right)\right] \right\}  \label{eq:twostage_objfunc}\\
\textrm{s. t.}\quad\sum_{t \in T^{d}_f}u_{f,t} &= 1, \forall f \in F, \label{eq:twostage_1b}\\
\sum_{t \in T^{a}_f}v_{f,t} &= 1, \forall f \in F, \label{eq:twostage_1c}\\
g_{f'} + a_{f'} - s_{f'} &\leq a_{f},\forall (f, f') \in \mathcal{C},\label{eq:twostage_1d}
\end{alignat}
\label{eq:twostage_a}
\end{subequations}

\noindent
where $Q\left(u,\xi^{(g)}\right)$ is:

\begin{subequations}
\begin{alignat}{2}
\min_{y^{(g)}} \quad & \sum_{z \in Z}\sum_{t \in T} C_g y^{(g)}_{z,t}\left(\xi^{(g)}\right) \label{eq:twostage_21a}\\
\textrm{s. t.}\sum_{f:z^{g}_{f} = z} u_{f,t} &\leq D_{z,t}\left(\xi^{(g)}\right) + y^{(g)}_{z,t} \left(\xi^{(g)} \right), \label{eq:twostage_21b}\\ &\qquad\qquad\;\forall \left(z,t,\xi^{(g)}\right) \in Z \Cross T \Cross \Xi^{(g)}, \nonumber \\
y^{(g)}_{z,t}\left(\xi^{(g)}\right) &\geq 0, \qquad\forall \left(z,t,\xi^{(g)}\right) \in Z \Cross T \Cross \Xi^{(g)}, \label{eq:twostage_21c}
\end{alignat}
\label{eq:twostage_2a}
\end{subequations}

\noindent
and $Q\left(v,\xi^{(a)}\right)$ is analogous to $Q\left(u,\xi^{(g)}\right)$ but for the arrivals end. Note that we use different ambiguity set sizes between arrival and departure capacity distributions. 

The objective function \eqref{eq:twostage_objfunc} includes two inner maximization problems, which seek the worst-case distribution within each ambiguity set that maximizes the expected second stage cost. Constraints \eqref{eq:twostage_1b}-\eqref{eq:twostage_1d} are the first stage assignment and coupling constraints inherited from the standard, deterministic MAGHP.  
The first constraint in the second stage minimization problem ensures that, even if the airport capacity is reduced, the number of departing (or arriving) flights at time $t$ does not exceed the realized airport capacity, plus the total number of extra flights allowed to depart (or arrive) at time $t$, across all airports. When there is a drop in airport capacity, the two-stage model will optimally adjust delay allocations based on the weighting between airborne and ground delays (typically airborne delays are more expensive \cite{delay_cost}). 


\subsection{Scenario reduction and \textsc{dr-MAGHP} reformulation} \label{ssec:scenario_reduce}

Recall that we are making decisions across a prediction horizon of 12 hours, subdivided into 48 time periods. From the airport capacity distribution prediction models (Section \ref{ssec:airport_cap_dis_pred}), we are given new predicted distributions of the arrival and departure capacities (and hence, associated Wasserstein ambiguity sets) at each time, across all airports.


This approach results in a substantial increase in the number of scenarios, leading to poor scalability and computational intractability. 
To reduce the number of scenarios, we focus on clustering time periods together with similar capacity distributions. 
We utilize the Wasserstein distance to quantify the pairwise similarity between capacity distributions of two consecutive time periods. Highly similar time periods are grouped together, whereas a distinct time period (marked by a significant change in predicted airport capacity distributions) is demarcated when the pairwise similarity exceeds a predefined threshold. The representative capacity distribution for a group is its average (centroid) capacity distribution. 


After scenario reduction, we reformulate the \textsc{dr-MAGHP} in \eqref{eq:twostage_objfunc}-\eqref{eq:twostage_1d} by converting the inner second-stage maximization problem into a minimization problem, resulting in a semi-infinite program. We then apply discretization techniques to address the continuous supports of this semi-infinite program, resulting in a deterministic equivalent form for the \textsc{dr-MAGHP}. Due to page limitations, we omit the technical details and refer readers to \cite{drGHP}. 

\section{Experimental results and discussion}  

\subsection{Capacity distribution prediction}

\subsubsection{Data description}

We obtained airport throughput data from the US Department of Transportation's Bureau of Transportation Statistics (BTS), and weather data from the US National Oceanic and Atmospheric Administration's High-Resolution Rapid Refresh (HRRR) database. BTS provides detailed information for each flight, including the scheduled departure and arrival times, actual departure and arrival times, delay duration, cancellation status, among others. Using these data points and the procedure described in Section \ref{sssec:derive_actual_capacities}, we estimate the actual airport capacity at each airport of interest. In this study, we examine the FAA Core 30 airports, but both the prediction and optimization frameworks can be easily extended to a larger network of airports. 
HRRR provides weather data on a $3$ km $\times$ $3$ km grid covering all 50 US states, with a forecast horizon of up to 23 hours from the current hour. We collect data for the entire year of 2019. Each day is divided into 96 quarter-hour intervals, resulting in 35,040 time periods in total.
We use 2019 weather data from January 1-December 30 to train the capacity distribution prediction models. 
We reserve December 31, 2019 as the day across which we will forecast airport capacity distributions that will serve as inputs for the \textsc{dr-MAGHP} model. 

\subsubsection{Prediction experimental setup and results}

We designate December 31, 2019 as the test set to evaluate model performance. For training-validation sets, we randomly split the remaining data into a standard ratio of $80$:$20$. We then normalize all numeric features on the training set using min-max  normalization. 
We use the same normalization procedure for the test set as well. 
We give the 3-layer MLP architecture details in Table \ref{tab:MLP_setting}; hyperparameter turning is done via grid search, with a resultant learning rate of 0.0001, 300 epochs, and a batch size of 16.


We use three metrics---\emph{Root Mean Squared Error (RMSE)}, \emph{Coverage Rate (CR)}, and \emph{Average Confidence Interval Length (ACIL)}---to evaluate model performance. Briefly, RMSE measures model performance with respect to a single, likeliest predicted capacity, ignoring the rest of the predicted distribution. In contrast, CR and ACIL incorporates differences between the predicted and actual capacity distributions through the use of confidence intervals. The deep learning model attempts to predict airport capacity distributions that minimize the RMSE, increase the CR, and reduce the spread of the ACIL. 

\subsubsection{Capacity distribution prediction results}
Table \ref{tab:prediction_results} summarizes the model performance across the test set for the US Core 30 airports on December 31, 2019. A salient future research direction would be to tune specific ambiguity sets within the \textsc{dr-MAGHP} in response to the prediction performance at individual airports. Furthermore, these performance metrics can guide future prediction model enhancements (e.g., using refined prediction models for airports with large RMSEs such as IAD and CLT, or low CR with a large ACIL).

\begin{figure}[htbp]
    \centering
    \includegraphics[width=0.4\textwidth]{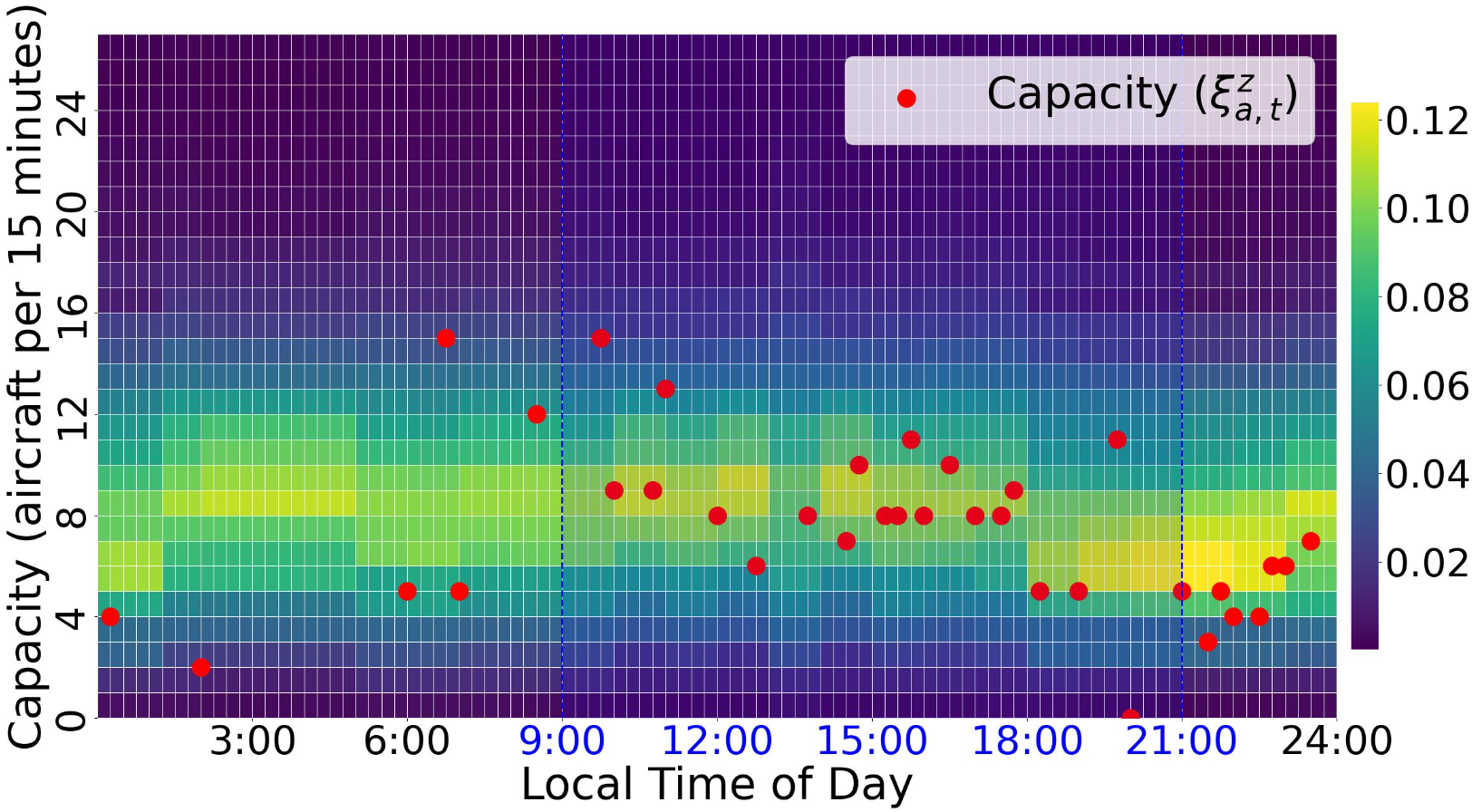}
    \caption{Predicted departure capacity distributions for LAX on December 31, 2019. Blue area from 9:00-21:00 is the 12-hour solution window of the \textsc{dr-MAGHP}.}
    \label{fig:lax_dep_heatmap_1231}
\end{figure}

We select Los Angeles International (LAX) Airport as a case study to visualize the prediction results. 
Figure~\ref{fig:lax_dep_heatmap_1231} plots the actual estimated departure capacity (red dots) and the predicted departure capacity distribution for LAX across December 31, 2019. 
At 18:00 local time, we observe a rapid decline in the predicted capacity, which aligns with the estimated actual capacity: During this time, LAX experienced a sudden change in the wind direction---this most likely triggered a change in the airport runway configuration. 
Furthermore, the predicted capacity distributions appear to be relatively stable, which reinforces the validity of the scenario reduction procedure described in Section \ref{ssec:scenario_reduce}. 

\begin{table}[htp]
    \centering
    \resizebox{0.45\textwidth}{!}{
    \begin{tabular}{c|ccc|ccc}
    \hline 
        \textbf{Airport} & \multicolumn{3}{c|}{\textbf{Arrival}} & \multicolumn{3}{c}{\textbf{Departure}} \\
        & \textbf{RMSE} & \textbf{CR (\%)} & \textbf{ACIL ($\mu\pm\sigma$)} & \textbf{RMSE} & \textbf{CR (\%)} & \textbf{ACIL ($\mu\pm\sigma$)} \\ \hline
        ATL & 5.00  & 100 & 25.95 $\pm$ 0.92 & 5.43  & 100 & 26.46 $\pm$ 0.84 \\
        BOS & 3.21  & 91  & 10.18 $\pm$ 0.39 & 2.42  & 100 & 11.83 $\pm$ 0.37 \\
        BWI & 1.78  & 100 & 7.50 $\pm$ 0.50 & 1.73  & 100 & 9.50 $\pm$ 0.50 \\
        CLT & 6.41  & 100 & 22.53 $\pm$ 0.75 & 11.61 & 93  & 24.00 $\pm$ 0.00 \\
        DCA & 2.72  & 100 & 10.20 $\pm$ 0.40 & 2.65  & 100 & 11.14 $\pm$ 0.35 \\
        DEN & 7.39  & 100 & 23.39 $\pm$ 0.89 & 9.09  & 86  & 26.86 $\pm$ 1.12 \\
        DFW & 8.87  & 96  & 25.91 $\pm$ 1.56 & 6.82  & 100 & 28.47 $\pm$ 0.70 \\
        DTW & 6.55  & 100 & 19.85 $\pm$ 1.66 & 8.85  & 91  & 22.36 $\pm$ 0.88 \\
        EWR & 2.59  & 94  & 10.35 $\pm$ 0.48 & 3.33  & 95  & 12.00 $\pm$ 0.00 \\
        FLL & 3.69  & 91  & 8.36 $\pm$ 1.23 & 3.12  & 100 & 8.42 $\pm$ 0.64 \\
        HNL & 2.05  & 100 & 8.00 $\pm$ 0.94 & 3.11  & 100 & 9.00 $\pm$ 0.00 \\
        IAD & 11.12 & 100 & 18.33 $\pm$ 0.47 & 7.31  & 100 & 21.20 $\pm$ 2.71 \\
        IAH & 6.77  & 100 & 19.87 $\pm$ 0.81 & 9.59  & 100 & 22.73 $\pm$ 0.45 \\
        JFK & 2.38  & 100 & 9.00 $\pm$ 0.00 & 3.44  & 100 & 13.57 $\pm$ 1.68 \\
        LAS & 2.75  & 100 & 10.57 $\pm$ 0.62 & 2.69  & 100 & 11.14 $\pm$ 0.64 \\
        LAX & 2.77  & 100 & 15.05 $\pm$ 0.86 & 3.42  & 95  & 15.77 $\pm$ 0.42 \\
        LGA & 3.21  & 100 & 11.33 $\pm$ 0.47 & 2.87  & 100 & 12.29 $\pm$ 0.45 \\
        MCO & 3.13  & 86  & 9.29 $\pm$ 0.96 & 2.06  & 100 & 10.58 $\pm$ 0.64 \\
        MDW & 1.00  & 100 & 6.25 $\pm$ 0.43 & 2.68  & 88  & 7.88 $\pm$ 0.32 \\
        MEM & 0.71  & 100 & 7.00 $\pm$ 0.00 & 2.00  & 100 & 8.00 $\pm$ 0.00 \\
        MIA & 3.52  & 92  & 10.25 $\pm$ 0.83 & 1.47  & 100 & 12.33 $\pm$ 1.60 \\
        MSP & 4.11  & 100 & 17.44 $\pm$ 0.50 & 5.59  & 100 & 21.00 $\pm$ 0.00 \\
        ORD & 3.73  & 100 & 25.48 $\pm$ 0.50 & 6.46  & 97  & 28.34 $\pm$ 1.59 \\
        PHL & 2.26  & 100 & 13.09 $\pm$ 0.51 & 5.67  & 80  & 13.70 $\pm$ 0.90 \\
        PHX & 3.70  & 100 & 14.84 $\pm$ 1.09 & 4.09  & 94  & 14.59 $\pm$ 0.49 \\
        SAN & 1.80  & 100 & 8.11 $\pm$ 0.99 & 2.85  & 100 & 11.80 $\pm$ 0.40 \\
        SEA & 2.59  & 100 & 13.19 $\pm$ 0.39 & 2.35  & 100 & 14.94 $\pm$ 0.23 \\
        SFO & 3.38  & 93  & 11.93 $\pm$ 0.59 & 3.48  & 93  & 13.14 $\pm$ 0.35 \\
        SLC & 2.81  & 100 & 15.78 $\pm$ 1.47 & 3.35  & 100 & 20.75 $\pm$ 0.83 \\
        TPA & 1.63  & 100 & 6.00 $\pm$ 1.15 & 2.52  & 100 & 10.00 $\pm$ 1.63 \\
        \hline 
    \end{tabular}}
    \caption{Capacity distribution prediction performance for FAA 30 core airports on Dec 31, 2019. RMSE and ACIL are given in aircraft per 15 minutes.} 
    \label{tab:prediction_results}
\end{table}

\subsection{\textsc{dr-MAGHP} optimization}
\subsubsection{Experiment setup}
To ensure consistency with the capacity distribution prediction procedure, we also use BTS data from December 31, 2019 to construct flight schedules. 
The solve horizon for the \textsc{dr-MAGHP} extends over twelve hours (09:00:00 to 21:00:00 local time), and is divided into forty-eight time periods ($|T| = 48$), each representing fifteen minutes. We optimize over 10,869 flights associated with the FAA Core 30 airports ($|Z| = 30$). 
We set realistic minimum turnaround times (45 minutes), and follow standard MAGHP setups such as an additional time period to absorb excess flights. 

\subsubsection{Results and comparisons}
We define the outputs from the upstream capacity distribution prediction procedure as predicted capacity distributions. To test the performance of the \textsc{DR-MAGHP} when predicted capacity distributions are not true capacity distributions, we define the realized capacity distributions as testing distributions. Testing distributions of interest are the cases where realized capacity distributions are shifted to the left, i.e., the realized probabilistic capacities are lower than anticipated. We also develop the stochastic programming based MAGHP (\textsc{sp-MAGHP}), which is the MAGHP model generating ground holding policies based on predicted capacity distributions, to compare its delay costs with those of the \textsc{DR-MAGHP} when predicted capacity distributions are not accurate. We use terms \emph{in-sample performance} ($\phi_{IS}$) and \emph{out-of-sample performance} ($\phi_{OS}$) to refer to the costs (i.e., optimal value in expectation) of \textsc{sp-MAGHP} and \textsc{dr-MAGHP} when evaluated on predicted distributions and testing distributions respectively.



\begin{figure}
    \centering
    \includegraphics[width=0.3\textwidth]{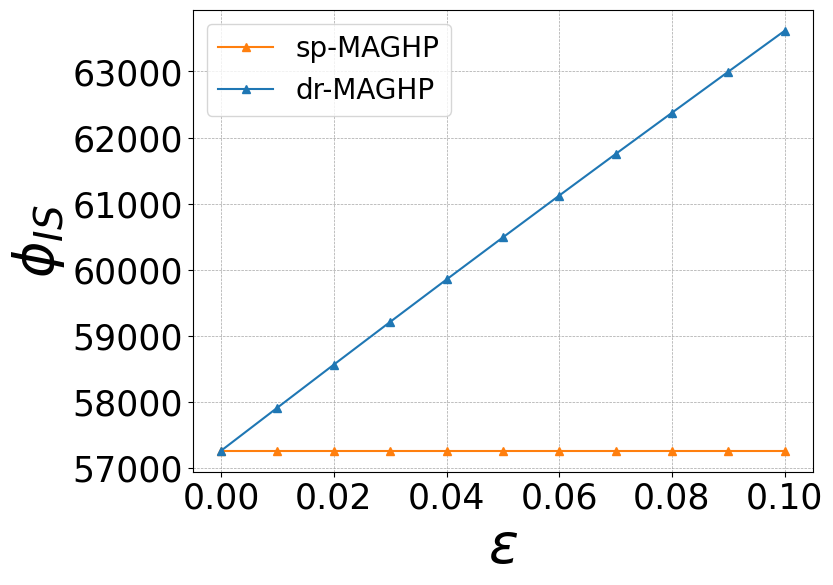}
    \caption{In-sample performance of \textsc{sp-MAGHP} and \textsc{dr-MAGHP} based on predicted capacity distributions}
    \label{fig:training_results}
\end{figure}

Figure \ref{fig:training_results} plots the in-sample performance of \textsc{dr-MAGHP} compared to \textsc{sp-MAGHP} as the size of the ambiguity sets grow (i.e., increasing radii \(\epsilon^{(a)}, \epsilon^{(g)}\)). When \(\epsilon^{(a)} = \epsilon^{(g)} = 0\), the ambiguity sets are singletons containing only the predicted distributions. Hence, the in-sample performance of \textsc{sp-MAGHP} is identical to that of \textsc{dr-MAGHP} when the ambiguity set contains only the predicted distribution, as expected. As the ambiguity sets grow, \textsc{dr-MAGHP} tends to generate more conservative ground holding policies, which increases delay costs when compared to \textsc{sp-MAGHP}.

\begin{algorithm}
\caption{Capacity Resampling Algorithm}
\begin{algorithmic}[1]
    \STATE \textbf{Input}: Predicted PMF $\{\widehat{p}^{(1)},\widehat{p}^{(2)},\dots, \widehat{p}^{(|Z|)}\}$ of each airport's capacity $\{\widehat{\xi}^{(1)},\widehat{\xi}^{(2)},\dots, \widehat{\xi}^{(|Z|)}\}$ with mean values $\{\mu^{(1)},\mu^{(2)},\dots,\mu^{|Z|}\}$, capacity reduction level $r$ for all airports and a maximum variability rate of probability $\delta$.
        \FOR{each airport $z$}
            \STATE Target mean $\mu^{*} = \widehat{\mu}^{(z)}  \left(1-r\right)$
            \STATE Update weights for each support $p_{(z)}$:
            \begin{equation}
            \begin{aligned}
                p^{(z)} = \arg\min_{p} \sum_{i=1}^{|\widehat{\xi}^{(z)}|}p_{i}\widehat{\xi}^{(z)}_{i} & \\
                \textrm{s.t.} \quad \sum_{i=1}^{|\widehat{\xi}^{(z)}|} p_{i}\widehat{\xi}^{(z)}_{i} &\geq \mu^{*},\quad\sum_{i=1}^{|\widehat{\xi}^{(z)}|}p_{i} = 1\\
                p_{i} - \widehat{p}^{(z)}_{i} \leq \delta\widehat{p}^{(z)}_{i},\quad-\delta\widehat{p}^{(z)}_{i} &\leq p_{i} - \widehat{p}^{(z)}_{i},\quad \forall i \in N
                \label{eq:lp_capreduc}
            \end{aligned}
            \end{equation}
        \ENDFOR
\STATE Draw i.i.d samples $\Bar{\xi}^{(z)} \sim p^{(z)}, z = 1,2,\dots,|Z|$ \\
\STATE \textbf{Output}: Reduced testing capacities for all airports $\{\Bar{\xi}^{(1)},\Bar{\xi}^{(2)},\dots,\Bar{\xi}^{(|Z|)}\}$
\end{algorithmic}
\label{Alg:alg1}
\end{algorithm}

\begin{figure}
\centering
    \subfloat[\label{subfig:sp_dr_0.1_increase}\textcolor{white}{Percent increase in flights subject}]{\includegraphics[width=4.42cm]{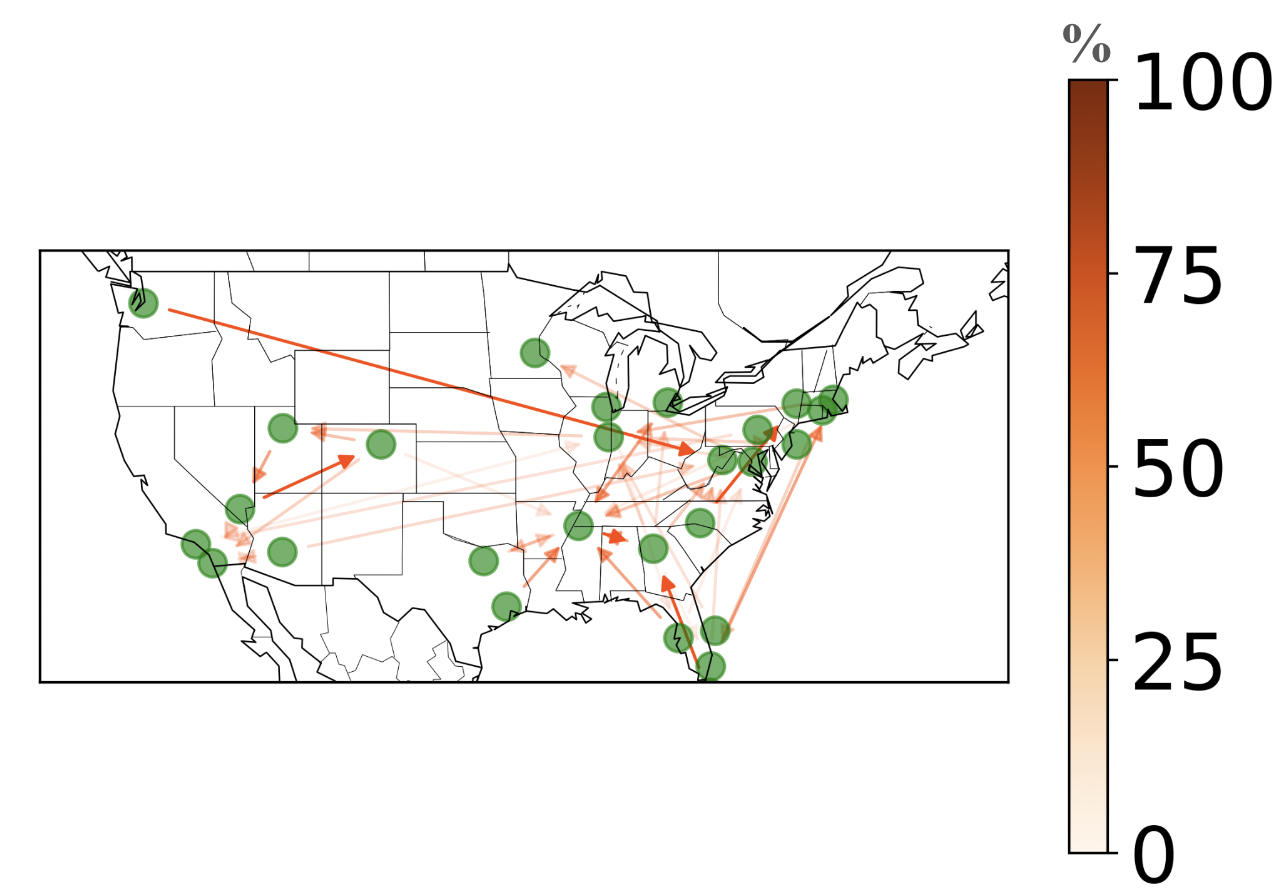}}\hfill
    \subfloat[\label{subfig:sp_dr_0.1_decrease}\textcolor{white}{Percent decrease in flights subject}]{\includegraphics[width=4.42cm]{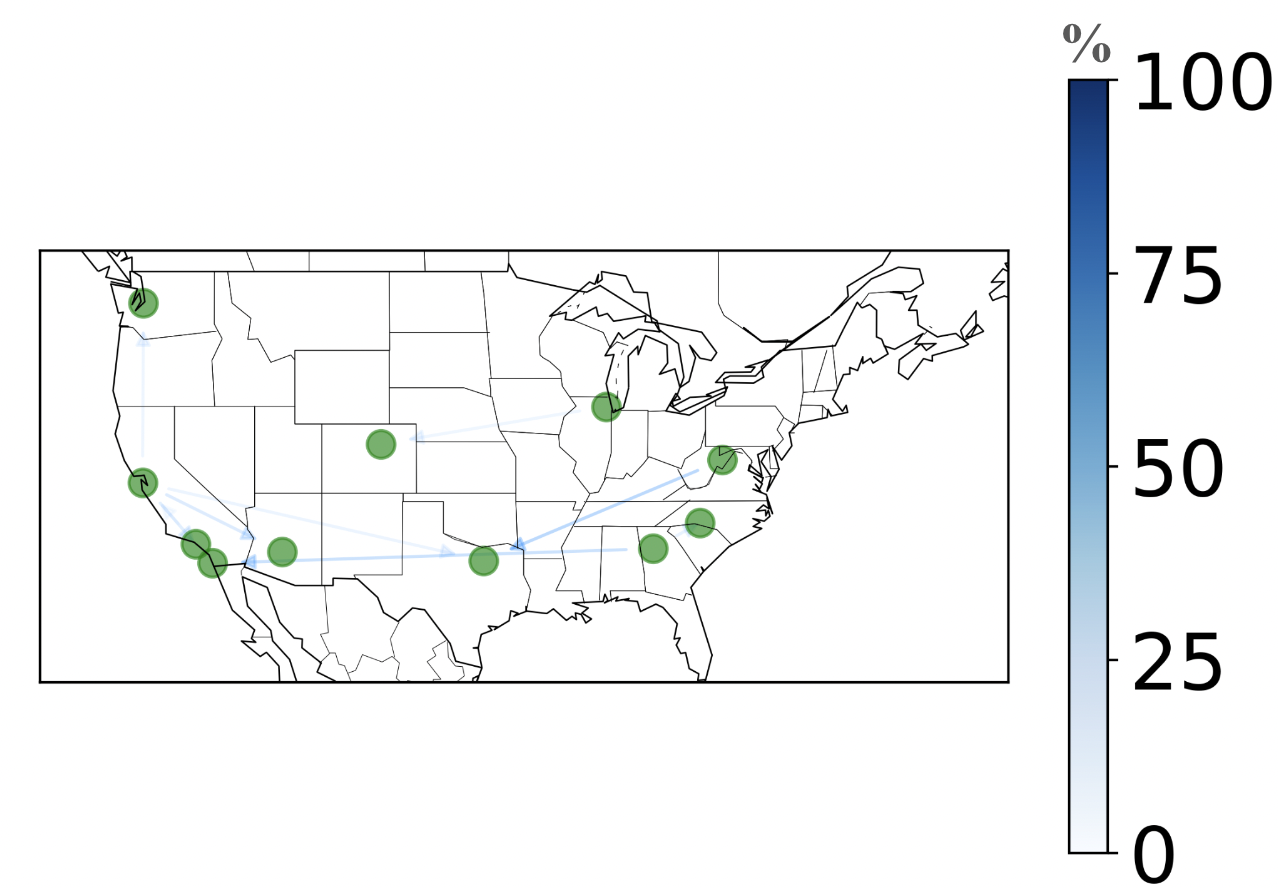}}\vfill
    \subfloat[\label{subfig:sp_dr_0.5_increase}\textcolor{white}{Percent increase in flights subject}] {\includegraphics[width=4.42cm]{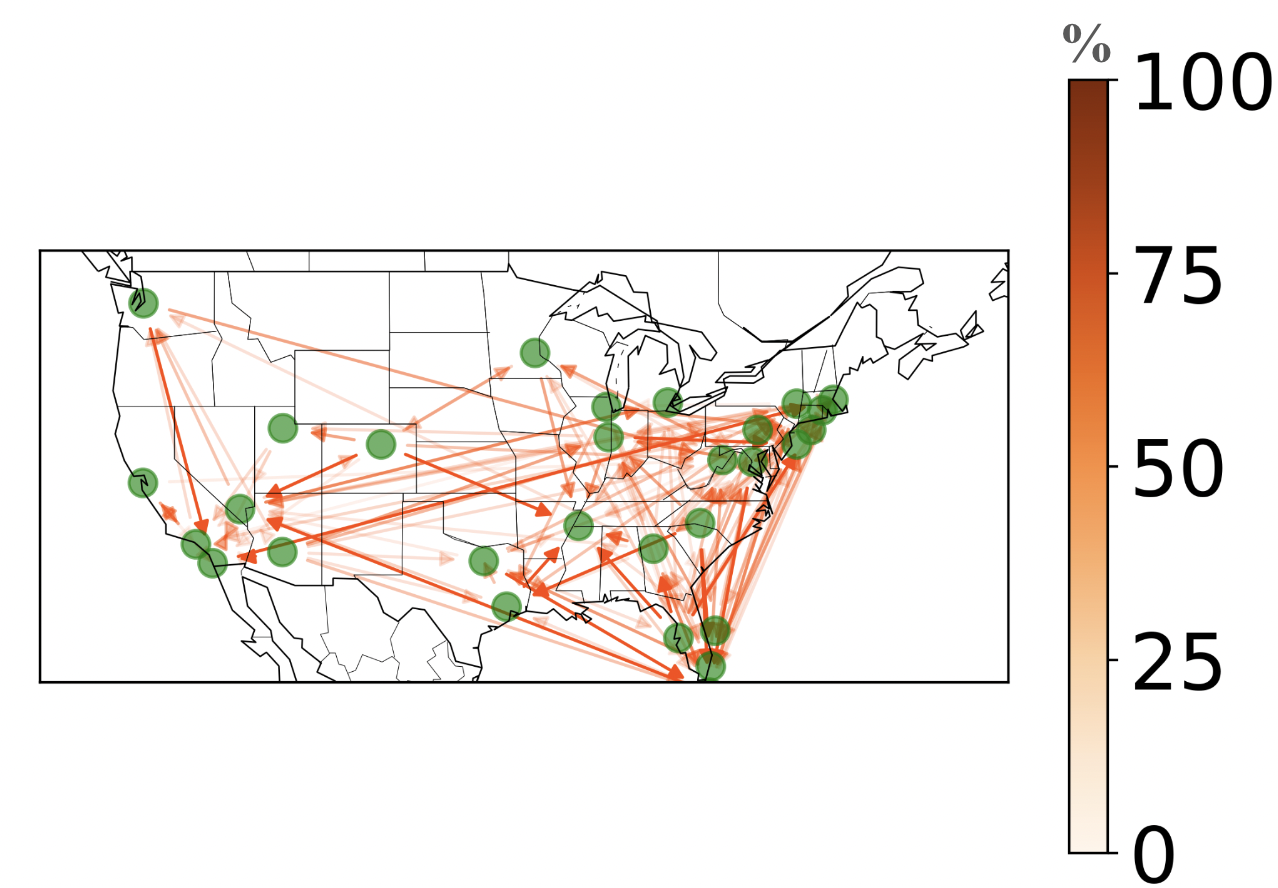}}\hfill
    \subfloat[\label{subfig:sp_dr_0.5_decrease}\textcolor{white}{Percent decrease in flights subject}]{\includegraphics[width=4.42cm]{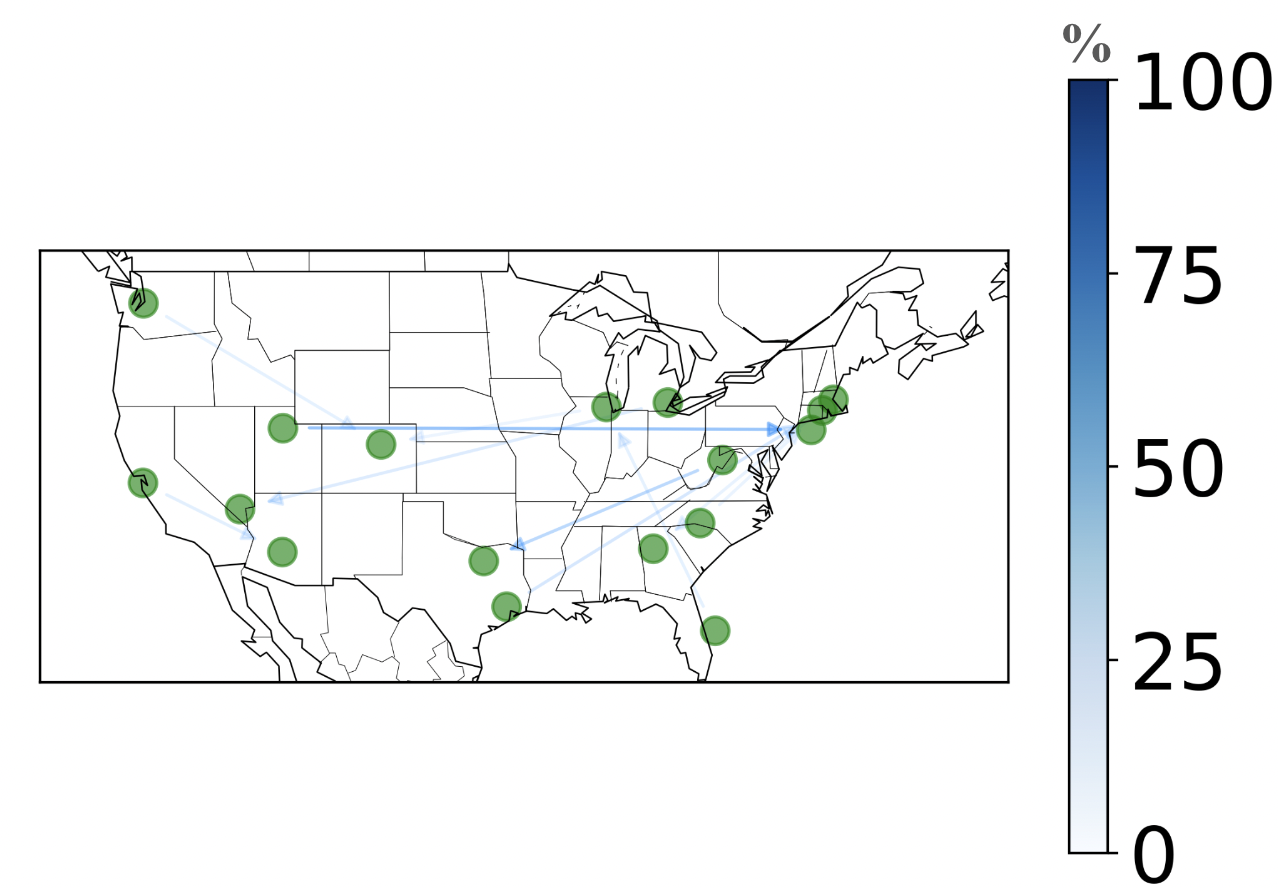}}\vfill
\caption{Percent \emph{increase} in delayed flights under \textsc{dr-MAGHP} versus \textsc{sp-MAGHP} with (a) $\epsilon = 0.1$, (c) $\epsilon = 0.5$; Percent \emph{decrease} in delayed flights under \textsc{dr-MAGHP} versus \textsc{sp-MAGHP} with (b) $\epsilon = 0.1$, (d) $\epsilon = 0.5$.}\label{fig:sp_dr_training}
\end{figure}

We provide a more operationally interpretable plot in Figure \ref{fig:sp_dr_training}, which depicts how GDP-induced delays change between the \textsc{sp-MAGHP} and \textsc{dr-MAGHP} solutions. Each directional edge represents the percentage of flights delayed across all direct flights from one airport node to another. Comparing Figures \ref{subfig:sp_dr_0.1_increase} and \ref{subfig:sp_dr_0.5_increase} reveals that larger ambiguity sets lead to more assigned delays across the board. 
We emphasize that these comparisons and interpretations must be made on a day-by-day basis, just as unique GDP policies must be developed for each encountered NAS scenario.

\subsubsection{Sensitivity to distribution shifts}
To evaluate how variations in the discrepancies between predicted and realized capacity distributions influence delay costs, and to determine the effectiveness of \textsc{dr-MAGHP} in mitigating these effects, we perform a sensitivity analysis. This analysis generates testing distributions at various levels of capacity reductions to compare the out-of-sample performance of \textsc{sp-MAGHP} and \textsc{dr-MAGHP}.




To sample from reduced capacity distributions at various reduction levels, we introduce a linear program in \eqref{eq:lp_capreduc} that performs valid adjustments to the PMFs 
to minimize the deviation between the current mean value and the targeted (reduced) mean value. 
We also introduce a parameter $\delta$ for the maximum variability rate to ensure a more uniform distribution of probability mass. 
The details of this procedure are outlined in Algorithm \ref{Alg:alg1}. We then draw 100 samples from the reduced capacity distributions, and evaluate the out-of-sample performance $\phi_{OS}(\Bar{x}) = \sum_{i=1}^{|\Bar{\xi}|}\widehat{\phi_{OS}}(\Bar{x},\Bar{\xi_{i}}) / |\Bar{\xi}| $, where $\widehat{\phi_{OS}}(\Bar{x},\Bar{\xi_{i}})$ is the objective function value of the derived ground holding policy $\Bar{x}$ with reduced capacity sample $\Bar{\xi_{i}}$.


We examine erroneous predictions where realized capacities are reduced between 10-100\% (note that a GDP with 0 AAR is essentially a Ground Stop). When the realized capacities are only slightly lower than forecast (i.e., 10-20\% reductions), only a small ambiguity set (i.e., small $\epsilon$) is needed to modestly outperform \textsc{sp-MAGHP}. As the realized capacities move farther away from predicted ones, larger ambiguity sets are needed. It is important to note that the out-of-sample performances of \textsc{sp-MAGHP} and \textsc{dr-MAGHP}, when radii are set to zero, might not align, even if their in-sample performances appear identical under the same conditions. This discrepancy arises because the specific delay assignment policies implemented by \textsc{sp-MAGHP} and \textsc{dr-MAGHP} based on predicted distributions may vary. These differences can lead to identical in-sample performance, yet diverge in out-of-sample performance.


Table \ref{tab:sensitivity_analysis_table} gives the out-of-sample performances (e.g., delay assignment costs for \textsc{sp-MAGHP} and \textsc{dr-MAGHP}, respectively) for different realized capacity distributions with reduced capacities. The best ambiguity set size $\epsilon^*$ is also given. Although the cost reductions may seem modest when utilizing \textsc{dr-MAGHP}, we note that this was across an arbitrary day (December 31, 2019) in the NAS. A full analysis of the advantages gained by robustifying ground holding decisions with ML-driven airport capacity prediction inputs will require performing this experiment across a longer time horizon. The goal of this paper is to provide the foundations for robust, learning-driven optimization approaches to strategic ATM considering uncertain airport capacities. Should users trust predictive models, then the robustness aspects can (and should, to avoid overly-conservative solutions) be dialed down. On the other hand, should users \emph{not} trust prediction model outputs (e.g., current conditions are particularly unstable, or the prediction horizon is long), distributionally robust approaches can be considered.


\begin{table}
\centering
\begin{tabular}{c|c c c c}
\hline
\text{Reduction} & $\phi_{OS}^{sp}$ & $\phi_{OS}^{dr}$ & $\epsilon^{*}$ & \% dec.  \\
\hline
$\textbf{10\%}$  & 331469.45 &  331443.65 & 0.01 & 0.007\\
$\textbf{20\%}$  & 395454.05 & 395420.20 &  0.01 & 0.009\\
$\textbf{30\%}$  & 473078.15 & 473068.50 & 0.02 & 0.002\\
$\textbf{40\%}$  & 566420.20 & 566335.80 & 0.02 & 0.015\\
$\textbf{50\%}$  & 673100.10 & 672943.65 & 0.04 & 0.023\\
$\textbf{60\%}$  & 792804.80 & 792645.95 & 0.02 & 0.020\\
$\textbf{70\%}$  & 935837.35 & 935473.75 & 0.04 & 0.038\\
$\textbf{80\%}$  & 1105686.30 & 1104923.10 & 0.09 & 0.069\\
$\textbf{90\%}$  & 1278823.30 & 1277608.70 & 0.09 & 0.095 \\
$\textbf{100\%}$ & 1454765.00 & 1453077.00 & 0.09 & 0.116\\
\hline
\end{tabular}
\caption{Sensitivity analysis results.}
\label{tab:sensitivity_analysis_table}
\end{table}

\begin{figure}
    \centering
    \subfloat[\label{subfig:0.1}$10\% $ reduction in capacities.]{\includegraphics[width=4.4cm]{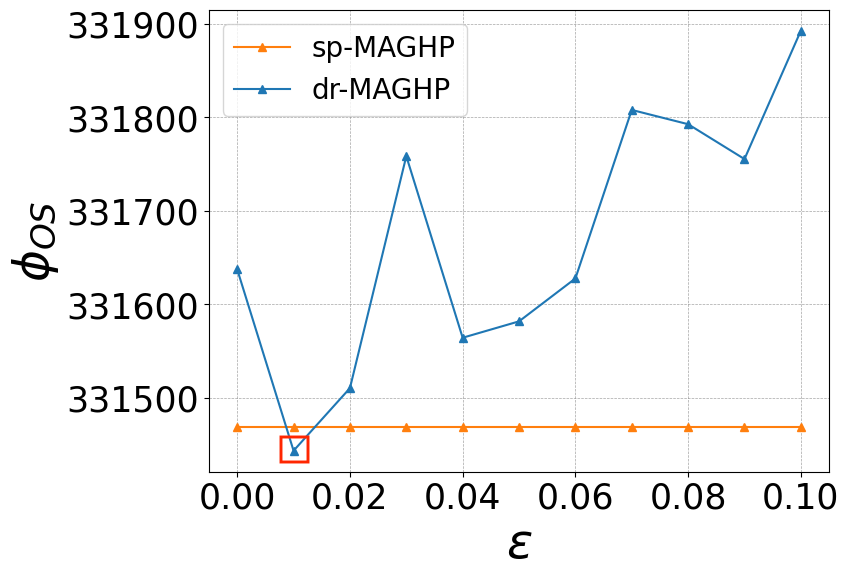}}\hfill
    \subfloat[\label{subfig:0.5}$50\%$ reduction in capacities.]{\includegraphics[width=4.4cm]{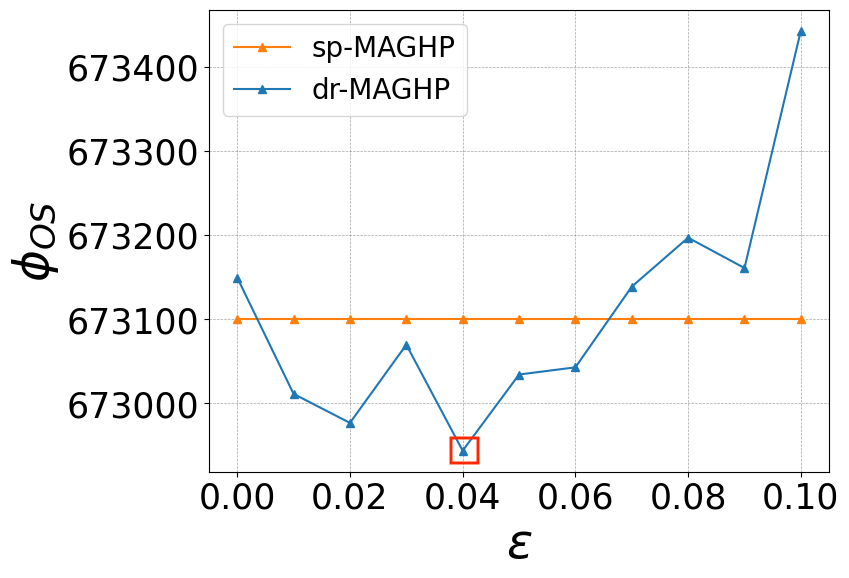}}\hfill
\caption{Out-of-sample performances of \textsc{sp-MAGHP} and \textsc{dr-MAGHPs} for two levels of capacity reductions.}\label{fig:sensitivity_analysis}
\end{figure}





\section{Concluding remarks}   

In this paper, we introduce a framework that merges an upstream airport capacity prediction model with a downstream distributionally robust ground delay program optimization model. The output of this learning-driven optimization framework is a policy for assigning ground and airborne delays to a set of flights, robustified against distributional uncertainties in the predicted capacities. We use a neural network approach (multilayer perceptron model) to produce distributional predictions of future airport capacities: The prediction is not a single capacity, but instead a set of different capacity scenarios, each with a probability of occurrence. These distributional predictions are fed into the constraints of \textsc{dr-MAGHP}, where we explore the impacts of ambiguity set sizes and prediction errors (in the form of distribution shifts). Future work includes expanding the prediction model inputs to include additional operational factors (e.g., runway configurations), improving the estimation method for deriving \enquote{true} airport capacities, and determining optimal ambiguity set sizes without exhaustively searching through a large range of $\epsilon$. Additionally, methodological improvements will be needed to consider non-independent arrival and departure processes both in the prediction and the optimization steps.

\bibliographystyle{IEEEtran} 
\small{
\bibliography{ICRAT2024_predictDRO/main.bib}
}



\end{document}